\documentclass[11pt]{article}
\usepackage{amsmath}
\usepackage{pifont}
\usepackage{amssymb}
\usepackage{amsfonts}
\usepackage{amscd}
\usepackage{epsfig}
\usepackage{amsmath,amstext,amsthm,amsbsy,amssymb}
\usepackage{graphicx}
\usepackage{epsfig}

\newcommand{\bc}{{\mathbb C}}
\newcommand{\br}{{\mathbb R}}
\newcommand{\bh}{{\mathbb H}}

\newcommand{\bp}{{\mathbb P}}
\newcommand{\mtx}[4]

\newtheorem{thm}{Theorem}[section]

\newtheorem{lem}[thm]{Lemma}
\newtheorem{prop}[thm]{Proposition}
\newtheorem{cor}[thm]{Corollary}

\newtheorem{dfn}{Definition}[section]

\newtheorem{example}{Example}[section]

\begin{document}

\title{J{\o}rgensen's Inequality and Collars in \\ $n$-dimensional Quaternionic
Hyperbolic Space}
\author{Wensheng Cao
\thanks{Supported by NSFs of China (No.10801107,No.10671004 ) and
NSF of Guangdong Province (No.8452902001000043)} \,   \\
Department of Mathematics,\\
 Wuyi University,
Jiangmen, \\
Guangdong 529020, P.R. China\\
e-mail: {\tt wenscao@yahoo.com.cn}\\
 \and
John R.\ Parker\\
Department of Mathematical Sciences,\\
Durham University, \\
Durham DH1 3LE, England\\
{email: {\tt j.r.parker@durham.ac.uk}}
 }
\date{}
\maketitle
\begin{abstract}
In this paper, we obtain analogues of J{\o}rgensen's
inequality for non-elementary groups of isometries of quaternionic
hyperbolic $n$-space generated by two elements, one of which is
loxodromic. Our result gives some improvement over earlier results
of Kim \cite{kim} and Markham \cite{mar}. These
results also apply to complex hyperbolic space and give improvements
on results of Jiang, Kamiya and Parker \cite{jkp}.

As applications, we use the  quaternionic  version of J{\o}rgensen's
inequalities to construct embedded collars about short, simple,
closed geodesics in quaternionic hyperbolic manifolds.  We show that
these canonical collars are disjoint from each other. Our results
give some improvement over earlier results of Markham and Parker and
answer an open question  posed in  \cite{mapa}.
\end{abstract}

\medskip


\section{Introduction}

J{\o}rgensen's inequality \cite{jor} gives a necessary
condition for a non-elementary two generator subgroup of
${\rm PSL}(2,\bc)$ to be discrete. As a quantitative version of
Margulis' lemma, this inequality has been generalised in many ways.
Viewing  ${\rm PSL}(2,\br)$, which is isomorphic to ${\rm PU}(1,1)$, as the 
holomorphic isometry group of complex hyperbolic 1-space, 
we can seek to generalise J{\o}rgensen's inequality to ${\rm PU}(n,1)$
for $n>1$, the holomorphic isometry group of higher dimensional 
complex hyperbolic space. Examples of this are the stable
basin theorem of Basmajian and Miner \cite{bm} (see also \cite{parsbt}) and
the complex hyperbolic J{\o}rgensen's inequality of Jiang, Kamiya and
Parker \cite{jkp}.

Kellerhals has generalised J{\o}rgensen's inequality to ${\rm PSp}(1,1)$.
This group is the isometry group of quaternionic hyperbolic 
1-space ${\bf H}_{\bh}^1$, which is the same as real 
hyperbolic 4-space ${\bf H}_{\br}^4$.
For more details of ${\rm PSp}(1,1)$, including a classification of the
elements, see \cite{cpw}.
It is interesting to seek generalisations of J{\o}rgensen's inequality
to ${\rm PSp}(n,1)$ for $n>1$, that is to higher dimensional quaternionic 
hyperbolic isometries. The first steps in this programme
were taken by Kim and Parker \cite{kimp} who gave a quaternionic
hyperbolic version of Basmajian and Miner's stable basin theorem.
Subsequently, Markham \cite{mar} and Kim \cite{kim}
independently gave versions of J{\o}rgensen's inequality for ${\rm PSp}(2,1)$.
Cao and Tan \cite{caog} obtained an analogue
of J{\o}rgensen's inequality for non-elementary groups of isometries
of quaternionic hyperbolic $n$-space generated by two elements, one of
which is elliptic.

In this paper we consider subgroups of ${\rm PSp}(n,1)$ with a loxodromic 
generator. Any loxodromic element $g$ of ${\rm PSp}(n,1)$ 
can be conjugated in ${\rm Sp}(n,1)$ to the form:
\begin{equation}\label{loxo}
{\rm diag}\Bigl(\lambda_1,\  \lambda_2,\  \cdots,\  \lambda_{n-1},
\lambda_n,\overline{\lambda}_n^{-1}\Bigr),
\end{equation}
where $\lambda_i\in\bh$ for $i=1,\,\ldots,\,n$  and $\overline{\lambda}_n^{-1}$ 
are right eigenvalues of $g$ with $|\lambda_i|=1$ for $i=1,\,\ldots,\,n-1$
and $|\lambda_n|>1$. We want to consider loxodromic maps that are close
to the identity. To make this precise, if $g\in {\rm Sp}(n,1)$ is a 
loxodromic map conjugate to (\ref{loxo}), we define the following conjugacy 
invariants:
\begin{equation}\label{mg} \delta(g)=\max\{|\lambda_i-1|: \ i=1,\cdots,n-1\},
\qquad
M_g=2\delta(g)+|\lambda_n-1|+|\overline{\lambda}_n^{-1}-1|
\end{equation}
Observe that $M_g>0$ and that the smaller $M_g$ is the closer $g$ is to
the identity. 
Note that $M_g$ is a natural generalisation of the invariant
$$
M_g=2|\lambda_1-1|+|\lambda_2-1|+|\overline{\lambda}_2^{-1}-1|.
$$ 
defined independently by Kim \cite{kim} and Markham \cite{mar} 
for ${\rm Sp}(2,1)$.

We consider groups generated by $g$ and $h$ that are close to each
other. To make this precise, we use the cross ratio of the fixed points of 
the two loxodromic maps $g$ and $hgh^{-1}$. We define the cross ratio in 
Section \ref{sec-background}. The statement of our main theorem is:

\begin{thm}\label{thmcp}
Let $g$ be a loxodromic element of ${\rm Sp}(n,1)$ with 
$M_g<1$ and with fixed points $u,\,v\in \partial {\bf H}_{\bh}^n$. 
Let $h$ be any other element of ${\rm Sp}(n,1)$. If
\begin{equation}\label{cond}
\bigl|[h(u),u,v, h(v)]\bigr|^{1/2}\bigl|[h(u),v,u, h(v)]\bigr|^{1/2}
<\frac{1-M_g}{M_g^2}
\end{equation}
then the group $\langle g,h \rangle$ is either elementary or not
discrete.
\end{thm}

We remark that this theorem is also valid for ${\rm SU}(n,1)$ and 
is stronger than both Theorems 4.1 and 4.2 of \cite{jkp}.
This theorem has some useful corollaries which we gather into a single result:

\begin{cor}\label{corcp}
Let $g$ be a loxodromic element of ${\rm Sp}(n,1)$ with 
$M_g<1$ and with fixed points $u,\,v\in \partial {\bf H}_{\bh}^n$. 
Let $h$ be any other element of ${\rm Sp}(n,1)$. 
Suppose that one of the following conditions holds:
\begin{eqnarray}
\bigl|[h(u),v,u,h(v)]\bigr|^{1/2} & < & \frac{1-M_g}{M_g}, \label{cond2a} \\
\bigl|[h(u),u,v,h(v)]\bigr|^{1/2} & < & \frac{1-M_g}{M_g}, \label{cond2b} \\
\bigl|[u,v,h(u),h(v)]\bigr|^{1/2} & < &  1-M_g, \label{cond-1.1} \\
\bigl|[h(u),u,v, h(v)]\bigr|+\bigl|[h(u),v,u, h(v)]\bigr| & < & 
\frac{2(1-M_g)}{M_g^2}. \label{condK2-improved}
\end{eqnarray}
Then the group $\langle g,h \rangle$ is either elementary or not
discrete.
\end{cor}

When $n=2$ the statement of Corollary \ref{corcp} with 
the conditions (\ref{cond2a}) and (\ref{cond2b}) was given independently 
by Kim, Theorem 3.1 of \cite{kim}, and Markham Theorem 1.1 of \cite{mar}
and for higher dimensions Cao gave these conditions in an earlier
preprint \cite{caoh}.
These results are a direct generalisation of Theorem 4.1 of \cite{jkp}.
They all follow from Theorem 2.4 of Markham and Parker \cite{mapa2}
and the observation (see the proof of Theorem \ref{nthm1.6} below) 
that for all ${\bf z}\in V_0$
\begin{equation}\label{eq-mp}
\Bigl|\bigl \langle g({\bf z}),{\bf z}\bigr\rangle\Bigr|
\le M_g \bigl|\langle {\bf z},{\bf u}\rangle\bigr|\,
\bigl|\langle {\bf z},{\bf v}\rangle\bigr|,
\end{equation}
which, in terms of the Cygan metric,  may be rewritten as equation
(10) of \cite{mapa2} with $d_g=|\lambda_n|^{1/2}$ and $m_g=M_g^{1/2}$.

The statement of Corollary \ref{corcp} with condition (\ref{condK2-improved})
is stronger than the corresponding results in dimension $n=2$ given 
by Kim and Markham.
Kim's criterion, Theorem 3.2 of \cite{kim}, is $M_g\leq \sqrt{2\sqrt{2}-1}-1$  
and 
$$
\bigl|[h(u),u,v, h(v)]\bigr|+\bigl|[h(u),v,u, h(v)]\bigr|<
\frac{2-2M_g-M_g^2+\sqrt{4-8M_g-8M_g^2-4M_g^3-M_g^4}}{2M_g^2}.
$$
Markham's criterion, Theorem 1.2 of \cite{mar}, is $M_g\leq \sqrt{2}-1$  and 
$$
\bigl|[h(u),u,v, h(v)]\bigr|+\bigl|[h(u),v,u, h(v)]\bigr|<
\frac{1-M_g+\sqrt{1-2M_g-M_g^2}}{M_g^2},
$$
which is a direct generalisation of Theorem 4.2 of \cite{jkp}.
It is easy to see that (when they are defined)
\begin{eqnarray*}
\frac{2(1-M_g)}{M_g^2} 
& > & \frac{1-M_g+\sqrt{1-2M_g-M_g^2}}{M_g^2} \\
& > & \frac{2-2M_g-M_g^2+\sqrt{4-8M_g-8M_g^2-4M_g^3-M_g^4}}{2M_g^2}.
\end{eqnarray*}
Therefore Kim and Markham's results follow from (\ref{condK2-improved}).

\medskip

Meyerhoff \cite{mey} used J{\o}rgensen's inequality to show that if
a simple closed geodesic in a hyperbolic 3-manifold is sufficiently
short, then there exists an embedded tubular neighbourhood of this
geodesic, called a {\sl collar}, whose width depends only on the length
(or the complex length) of the closed geodesic. Moreover, he showed
that these collars were disjoint from one another. In
\cite{kel1,kel2} Kellerhals generalised   Meyerhoff's results to
real hyperbolic 4-space and 5-space with the aid of some properties
of quaternions. Markham and Parker \cite{mapa} used the complex and
quaternionic hyperbolic J{\o}rgensen's inequality obtained in
\cite{jkp,mar}, to give analogues of Meyerhoff's (and Kellerhals')
results for short, simple, closed geodesics in 2-dimensional complex
and quaternionic hyperbolic manifolds.  They showed  that these
canonical collars are disjoint from each other and from canonical
cusps. For complex hyperbolic space, by using a  lemma of Zagier
they also gave an estimate based only on the length,  and left the
same question for the case of quaternionic space as an open
question.

Let $G$ be a discrete group of $n$-dimensional quaternionic hyperbolic
isometries. Let $g\in G$ be loxodromic with axis the geodesic
$\gamma$. The {\sl tube} $T_r(\gamma)$ of radius $r$ about $\gamma$ is the
collection of points a distance less than $r$ from $\gamma$. It is
clear that $g$ maps $T_r(\gamma)$ to itself.   The tube
$T_r(\gamma)$ is {\sl precisely invariant} under the subgroup 
$\langle g \rangle$ of $G$ if $h\bigl(T_r(\gamma)\bigr)$ is disjoint from 
$T_r(\gamma)$ for all $h\in G-\langle g \rangle$. 
If $T_r(\gamma)$ is precisely
invariant under $G$  then $C_r(\gamma')=T_r(\gamma)/\langle g \rangle$ 
is an embedded tubular neighbourhood of the simple closed
geodesic $\gamma'=\gamma/\langle g \rangle$. We call $C_r(\gamma')$
the {\sl collar} of width $r$ about $\gamma'$.

As applications of our quaternionic version
J{\o}rgensen's inequalities, we will give analogues of Markham and
Parker's results for short, simple, closed geodesics in
$n$-dimensional quaternionic hyperbolic manifolds. 

Given a loxodromic map $g$ with axis $\gamma$ and satisfying $M_g<\sqrt{3}-1$, 
we define a positive real number $r$ by:
\begin{equation}\label{eq-can-width}
\cosh(2r)=\frac{2(1-M_g)}{M_g^2}.
\end{equation}
Then we call the tube $T_r(\gamma)$ with $r$ given by (\ref{eq-can-width})
the {\sl canonical tube} about $\gamma$. If $\gamma'=\gamma/\langle g \rangle$
then we call the collar $C_r(\gamma')$ with $r$ given by (\ref{eq-can-width})
the {\sl canonical collar} about $\gamma'$.

\begin{thm}\label{nthm1.4} Let $G$ be a discrete, non-elementary, torsion-free
subgroup of ${\rm Sp}(n,1)$.  Let $g$ be a loxodromic element of $G$ with
axis the geodesic $\gamma$. Suppose that $M_g <\sqrt{3}-1$. Then 
the canonical tube $T_r(\gamma)$ whose width $r$ 
is given by (\ref{eq-can-width}) 
is precisely invariant under $\langle g\rangle$ in $G$.

In particular, the canonical collar $C_r(\gamma')$ of width $r$ about
$\gamma'=\gamma/\langle g \rangle$ is embedded 
in the manifold $\mathcal{M}={\bf H}^n_{\bh}/G$. 
\end{thm}

Furthermore, we have

\begin{thm}\label{nthm1.6} Let $\mathcal{M}$ denote a quaternionic hyperbolic
$n$-manifold. Then the canonical collars around distinct short, simple,
closed geodesics in  $\mathcal{M}$ are disjoint.
\end{thm}

By  controlling  the rotational  part of  loxodromic element,  we
obtain the radius of collars solely in terms of the length of
the corresponding simple closed geodesic as the following, which
answers the open problem posed in \cite{mapa}.

\begin{thm}\label{nthm1.7}
Let $N\geq 35$ be a positive integer. Let $G$ be a discrete,
torsion-free, non-elementary subgroup of ${\rm Sp}(n,1)$.  Let $g$ be a
loxodromic element of G with axis $\gamma$ having the form
(\ref{loxo}) and let $l=2\log|\lambda_n|$ be the length of the closed geodesic
$\gamma/\langle g\rangle$ and suppose that 
\begin{equation}\label{condi}
R_N=2\sqrt{\left(\cosh\frac{N^nl}{2}+1\right)
\left(\cosh\frac{N^nl}{2}-\cos\frac{2\pi}{N}\right)}+
2\sqrt{2\left(1-\cos\frac{2\pi}{N}\right)}<\sqrt{3}-1.
\end{equation}
Define the positive number $r$ by
$$
\cosh(2r)=\frac{2(1-R_N)}{R_N^2}.
$$
Then the tube $T_r(\gamma)$
is precisely invariant under $G$.
\end{thm}

\begin{cor}\label{cor-1.8}
 Let $G$ be a discrete,
torsion-free, non-elementary subgroup of ${\rm Sp}(2,1)$.  Let $g$ be a
loxodromic element of G with axis $\gamma$ having the form
(\ref{loxo}).  Suppose that $l=2\log|\lambda_2|<0.00017681$. Let
$r$ be a positive number defined by
$$
\cosh(2r)=\frac{2(1-R)}{R^2}
$$
where
\begin{equation}\label{condi6}
R=2\sqrt{\left(\cosh\frac{1849\,l}{2}+1\right)
\left(\cosh\frac{1849\,l}{2}-\cos\frac{2\pi}{43}\right)}+
2\sqrt{2\left(1-\cos\frac{2\pi}{43}\right)}.
\end{equation}
Then the tube $T_r(\gamma)$ is precisely invariant under $G$.
\end{cor}

The structure of the remainder of this paper is as follows. In
Section 2,  we give the necessary background material for
quaternionic hyperbolic space. 
Section 3 contains  the proof of Theorem \ref{thmcp}
and Corollary \ref{corcp}.
In Section 4, we use Theorem \ref{thmcp} to obtain
the proof of Theorems \ref{nthm1.4} and \ref{nthm1.6}.
In Section 5, we give an example to illustrate the idea behind 
Theorem \ref{nthm1.7}.  Using the adapted {\sl Pigeonhole Principle}
(cf. \cite{mey}),  we obtain the proof of Theorem \ref{nthm1.7} and
Corollary \ref{cor-1.8}.

\section{Background}\label{sec-background}

We begin with some background material on quaternionic
hyperbolic geometry. Much of this can be found in
\cite{chen,gold,kimp,mos}. 
Let $\bh^{n,1}$ be the quaternionic vector space of  quaternionic
dimension $n+1$  (so real dimension $4n+4$) with the quaternionic
Hermitian form
$$
\langle{\bf z},\,{\bf w}\rangle={\bf w}^*J{\bf z}=
\overline{w}_1z_1+\cdots+\overline{w}_{n-1}z_{n-1}
-(\overline{w}_{n}z_{n+1}+\overline{w}_{n+1}z_{n}),
$$
where ${\bf z}$ and ${\bf w}$ are the column vectors in $\bh^{n,1}$
with entries $z_1,\cdots,z_{n+1}$ and $w_1,\cdots,w_{n+1}$
respectively, $\cdot^*$ denotes quaternionic Hermitian transpose and
$J$ is the Hermitian matrix
$$
J=\left(
                  \begin{array}{ccc}
                    I_{n-1} & 0 & 0 \\
                    0 & 0 & -1 \\
                    0 & -1 & 0\\
                  \end{array}
                \right).
$$
We define a {\sl unitary quaternionic transformation} (or symplectic 
transformation) $g$ to be an automorphism of
$\bh^{n,1}$, that is, a linear bijection such that 
$\langle g({\bf z}),\,g({\bf w})\rangle=\langle{\bf z},\,{\bf w}\rangle$ for all
${\bf z}$ and ${\bf w}$ in $\bh^{n,1}$. We denote the group of all
unitary transformations by ${\rm Sp}(n,1)$.

Following Section 2 of \cite{chen}, let
\begin{eqnarray*}
V_0 & = & \Bigl\{{\bf z} \in  \bh^{n,1}-\{0\}:
\langle{\bf z},\,{\bf z}\rangle=0\Bigr\} \\
V_{-} &  = & \Bigl\{{\bf z} \in \bh^{n,1}:\langle{\bf z},\,{\bf
z}\rangle<0\Bigr\}.
\end{eqnarray*}
It is obvious that $V_0$ and $V_{-}$ are invariant under ${\rm Sp}(n,1)$.

We define an equivalence relation $\sim$ on $\bh^{n,1}$ by 
${\bf z}\sim{\bf w}$ if and only if there exists a non-zero 
quaternion $\lambda$ so that ${\bf w}={\bf z}\lambda$. Let $[{\bf z}]$ denote
the equivalence class of ${\bf z}$. Let 
$\bp:\bh^{n,1}-\{0\}\longrightarrow \bh\bp^n$ be the {\sl right projection} map
given by $\bp:{\bf z}\longmapsto [{\bf z}]$. If $z_{n+1}\neq 0$ then 
$\bp$ is given by
$$
\bp(z_1,\,\ldots,\,z_n,  z_{n+1})^t=(z_1z_{n+1}^{-1},\cdots,z_n
z_{n+1}^{-1})^t\in{\bh }^n.
$$
We also define 
\begin{equation}\label{infty}
\bp(0,\,\ldots,\,0,\,z_n,\,0)^t=\infty.
\end{equation}
Observe that
\begin{equation}\label{proj-form}
\langle {\bf z}\lambda,{\bf w}\mu\rangle
=\overline{\mu}{\bf w}^*J{\bf z}\lambda
=\overline{\mu}\langle {\bf z},{\bf w}\rangle\lambda.
\end{equation}

We define the Siegel domain model of {\sl quaternionic hyperbolic $n$ space} 
to be ${\bf H}_{\bh}^n=\bp(V_-)$ and its boundary to be 
$\partial {\bf H}_{\bh}^n=\bp(V_0)$.
It is clear that $\infty\in\partial {\bf H}_{\bh}^n$. Also for
all ${\bf z}\in V_-$ we have $z_{n+1}\neq 0$ and so $\bp$ is given by
the formula above. Likewise for all ${\bf z}\in V_0$, either
$z_{n+1}\neq 0$ or $\bp({\bf z})=\infty$.

As in Chapter 19 of \cite{mos}, the Bergman metric on 
${\bf H}_{\bh}^n$ is given by the distance formula
$$
\cosh^2\frac{\rho(z,w)}{2}=
\frac{\langle{\bf z},\,{\bf w}\rangle \langle{\bf w},\,{\bf z}\rangle}
{\langle{\bf z},\,{\bf z}\rangle \langle{\bf w},\,{\bf w}\rangle},
\ \ \mbox{where}\ \ z,w \in {\bf H}_{\bh}^n, \ \  
{\bf z}\in \bp^{-1}(z),{\bf w}\in \bp^{-1}(w).
$$
This expression is independent of the choice of ${\bf z}$ and ${\bf w}$.
Since ${\rm Sp}(n,1)$ preserves the form $\langle \cdot,\,\cdot\rangle$, 
it clearly preserves the right hand side of this expression. Therefore
$g\in{\rm Sp}(n,1)$ acts on ${\bf H}_\bh^n\cup\partial{\bf H}_\bh^n$ as follows:
$$
g(z)=\bp g\bp^{-1}(z).
$$
This formula is well defined provided
the action of ${\rm Sp}(n,1)$ is on the left
and the action of projection $\bp$ of ${\rm Sp}(n,1)$ is on the right.
It is clear that multiples of $g$ by a non-zero real number act in the
same way. Since elements of ${\rm Sp}(n,1)$ have determinant $\pm 1$
this real number can only be $\pm 1$.
Therefore we define 
${\rm PSp}(n,1)={\rm Sp}(n,1)/\{\pm I_{n+1}\}$.
All elements of ${\rm PSp}(n,1)$ are isometries of ${\bf H}_\bh^n$.
We often find it convenient to work with matrices in ${\rm Sp}(n,1)$
rather than projective mappings in ${\rm PSp}(n,1)$ and we will pass
between them without comment.

If $g\in {\rm Sp}(n,1)$,  by definition, $g$ preserves the Hermitian form.
Hence
$$
{\bf w}^*J{\bf z}=\langle{\bf z},\,{\bf w}\rangle= \langle g{\bf
z},\,g{\bf w}\rangle ={\bf w}^* g^*J g{\bf z}
$$
for all ${\bf z}$ and ${\bf w}$ in $\bh^{n,1}$. Letting ${\bf z}$
and ${\bf w}$ vary over a basis for $\bh^{n,1}$, we see that $J=
g^*J g$. From this we find $ g^{-1}=J^{-1} g^*J$. That is:
\begin{equation}\label{hform}
g^{-1}=\left(
  \begin{array}{ccc}
     A^*& -\theta^*& -\eta^* \\
    -\beta^* & \overline{d}& \overline{b}\\
    -\alpha^* & \overline{c}& \overline{a}\\
    \end{array}
\right) \ \ \mbox{for}\ \
 g=\left(
  \begin{array}{ccc}
     A& \alpha& \beta \\
    \eta & a& b\\
    \theta & c& d\\
    \end{array}
\right)\in {\rm Sp}(n,1).\end{equation}
 Using the identities $gg^{-1}=g^{-1}g=I_{n+1}$ we obtain:
\begin{eqnarray}
 \label{AA} AA^*-\alpha \beta^*-\beta \alpha^* & = & I_{n-1},\\
   -A\theta^*+\alpha \overline{d}+\beta \overline{c} & = &  0,\\
  -A\eta^*+\alpha \overline{b}+\beta \overline{a} & = &  0,\\
  \label{BB}-\eta\theta^*+a\overline{d}+b\overline{c} & = &  1,\\
 \label{gg}-\eta\eta^*+a\overline{b}+b\overline{a} & = & 0,\\
 \label{tt}-\theta\theta^*+c\overline{d}+d\overline{c} & = & 0,\\
 A^*A-\theta^* \eta-\eta^* \theta & = & I_{n-1},\\
 \label{rr} A^*\alpha-\theta^*a-\eta^*c & = &  0,\\
 A^*\beta-\theta^*b-\eta^* d & = &  0,\\
 \label{bba}-\beta^*\alpha+\overline{d}a+\overline{b}c & = &  1,\\
 \label{bb} -\beta^*\beta+\overline{d}b+\overline{b}d & = & 0,\\
 \label{aa}-\alpha^*\alpha+\overline{c}a+\overline{a}c & = & 0.
\end{eqnarray}
Following Chen and Greenberg \cite{chen}, we say that a non-trivial element 
$g$ of ${\rm Sp}(n,1)$ is:
\begin{itemize}
\item[(i)] {\sl elliptic} if it has a fixed point in  ${\bf H}_{\bh}^n$;
\item[(ii)] {\sl parabolic} if it has exactly one fixed point which lies in
$\partial {\bf H}_{\bh}^n$;
\item[(iii)] {\sl loxodromic} if it has exactly two fixed points which
lie in $\partial {\bf H}_{\bh}^n$.
\end{itemize}

A  subgroup $G$ of ${\rm Sp}(n,1)$ is called {\sl elementary} if it has
a finite orbit in ${\bf H}_\bh^n\cup\partial{\bf H}_\bh^n$. If all
of its orbits are infinite then $G$ is {\sl non-elementary}. 
In particular, $G$ is non-elementary if it contains two
non-elliptic elements of infinite order with distinct fixed points.

Let $o$ be the origin in $\bh^n$ and $\infty$ be as defined in (\ref{infty}).
Both these points lie on $\partial{\bf H}_\bh^n$.
In what follows we make fixed choices of points in $\bh^{n,1}$ that
are preimages of these points.
Namely
$$
(0,\ldots,0,0,1)^t\in\bp^{-1}(o)\subset V_0, \qquad
(0,\ldots,0,1,0)^t\in\bp^{-1}(\infty)\subset V_0.
$$
Define the stabilisers of the points to be:
$$
G_o=\{g\in {\rm Sp}(n,1): g(o)=o\},\quad 
G_{\infty}=\{g\in {\rm Sp}(n,1): g(\infty)=\infty\},\quad  
G_{o,\infty}=G_o\cap G_{\infty}.
$$
Note that if $g$ has the form (\ref{hform}) then if $g\in G_o$ we have
$b=0$ and if $g\in G_\infty$ we have $c=0$. 

Cross-ratios were generalised to complex hyperbolic space by
Kor\'anyi and Reimann \cite{km}. We will generalise this definition
of complex cross-ratio to the non commutative quaternion ring.

\begin{dfn}   The quaternionic
cross-ratio of four points  $z_1, z_2, w_1, w_2$ in
$\overline{\bf H}_{\bh}^n$  is defined as:
\begin{equation}
 [z_1, z_2, w_1, w_2]=
\langle{\bf w_1},\,{\bf z_1}\rangle 
\langle{\bf w_1},\,{\bf z_2}\rangle^{-1}
\langle{\bf w_2},\,{\bf z_2}\rangle 
\langle{\bf w_2},\,{\bf z_1}\rangle^{-1}, 
\end{equation}
where ${\bf z}_i=\in \bp^{-1}(z_i)$ and ${\bf w_i}\in \bp^{-1}(w_i)$
for $i=1,\,2$.
\end{dfn}
Using (\ref{proj-form}) we see that
\begin{eqnarray*}
 [z_1\lambda_1, z_2\lambda_2, w_1\mu_1, w_2\mu_2] & = &
\langle{\bf w_1}\mu_1,\,{\bf z_1}\lambda_1\rangle 
\langle{\bf w_1}\mu_1,\,{\bf z_2}\lambda_2\rangle^{-1}
\langle{\bf w_2}\mu_2,\,{\bf z_2}\lambda_2\rangle 
\langle{\bf w_2}\mu_2,\,{\bf z_1}\lambda_2\rangle^{-1} \\
& = & \overline{\lambda}_1\langle{\bf w_1},\,{\bf z_1}\rangle \mu_1
\mu_1^{-1}\langle{\bf w_1},\,{\bf z_2}\rangle^{-1}\overline{\lambda}_2^{-1}
\overline{\lambda}_2\langle{\bf w_2},\,{\bf z_2}\rangle \mu_2
\mu_2^{-1}\langle{\bf w_2},\,{\bf z_1}\rangle^{-1} \overline{\lambda}_1^{-1} \\
& = & \overline{\lambda}_1 [z_1, z_2, w_1, w_2] \overline{\lambda}_1^{-1}.
\end{eqnarray*}
The quaternionic cross-ratio $[z_1, z_2, w_1, w_2]$ depends on
the choice of ${\bf z}_1\in\bp^{-1}(z_1)$.
However, its absolute value
\begin{equation}
 \bigl|[z_1, z_2, w_1, w_2]\bigr|
=\frac{|\langle{\bf w_1},\,{\bf z_1}\rangle\langle{\bf w_2},\,{\bf z_2}\rangle|}
{| \langle{\bf w_1},\,{\bf z_2}\rangle\langle{\bf w_2},\,{\bf z_1}\rangle|}
\end{equation} 
is independent of the preimage of $z_i$ and $w_i$ in $\bh^{n,1}$. 
The following lemma is easy to prove.

\begin{lem}\label{lem-X-ratio}
Let  $o, \infty\in \partial {\bf H}_{\bh}^n$ stand for the images of
$(0,\cdots,0,1)^t$ and $(0,\cdots,0,1,0)^t\in \bh^{n,1}$ under  the
projection map $\bp$, respectively and let $h\in{\rm PSp}(n,1)$ be given by
(\ref{hform}). 
Then
\begin{eqnarray}
\bigl|[h(\infty),o,\infty, h(o)]\bigr| & = & |bc|, \label{ratio1} \\
\bigl|[h(\infty),\infty,o, h(o)]\bigr| & = & |ad|, \label{ratio2} \\
\bigl|[\infty,o, h(\infty), h(o)]\bigr| & = & \frac{|bc|}{|ad|}. \label{ratio3}
 \end{eqnarray}
\end{lem}

The following lemma is crucial for us to prove  Theorem \ref{thmcp}.

\begin{lem}
\label{lemcu}
 Let $h$ be as in (\ref{hform}). Then
\begin{eqnarray}
\label{ba}
|\beta^*\alpha| & \leq &  2|ad|^{1/2}|bc|^{1/2}, \\
\label{gt}
|\eta\theta^*| & \leq &  2|ad|^{1/2}|bc|^{1/2}, \\
\label{nn}
|ad|^{1/2} & \leq &  |bc|^{1/2}+1, \\
\label{nd} 
|bc|^{1/2} & \leq &  |ad|^{1/2}+1, \\
\label{cl1}
1 & \le & |ad|^{1/2}+|bc|^{1/2}.
\end{eqnarray}
\end{lem}

\noindent
{\bf Proof.} Using (\ref{bb}) and (\ref{aa}), we have
\begin{equation}\label{bbaa}
|\beta^*\alpha|^2 \leq |\beta^*\beta|\,|\alpha^*\alpha|
=2\Re(\overline{d}b)\,2\Re(\overline{c}a)\leq 4|ad||bc|.
\end{equation}
This gives (\ref{ba}). Similarly, using (\ref{gg}) and (\ref{tt}), we have
$$
|\eta\theta^*|^2\leq |\eta\eta^*|\,  |\theta\theta^*|
=2\Re(a\overline{b})\, 2\Re(c\overline{d})
\leq 4|ad||bc|.
$$
This gives (\ref{gt}).

Next, using (\ref{bba}) and  (\ref{bbaa}), we have 
\begin{eqnarray*}
4\Re(\overline{d}b)\,\Re(\overline{c}a)
& \ge & |\beta^*\alpha|^2 \\
& = & |\overline{d}a+\overline{b}c-1|^2 \\
& = & 1+|ad|^2+|bc|^2-2\Re(\overline{d}a)-2\Re(\overline{b}c)
+2\Re(\overline{d}a\overline{c}b).
\end{eqnarray*}
Thus
\begin{eqnarray*}
1+|ad|^2+|bc|^2 & \le & 2\Re(\overline{d}a)+2\Re(\overline{b}c)
-2\Re(\overline{d}a\overline{c}b)+4\Re(\overline{d}b)\,\Re(\overline{c}a) \\
& = & 2\Re(\overline{d}a)+2\Re(\overline{b}c)+2\Re(b\overline{d}c\overline{a}) 
\\
& \le & 2|ad|+2|bc|+2|ad|\,|bc|.
\end{eqnarray*}
We can rearrange this expression to obtain
$$
\bigl(1-|ad|-|bc|\bigr)^2 \le 4|ad|\,|bc|.
$$
Taking square roots gives
$$
-2|ad|^{1/2}|bc|^{1/2}\le 1-|ad|-|bc|\le 2|ad|^{1/2}|bc|^{1/2}.
$$
Rearranging gives
$$
\bigl(|ad|^{1/2}-|bc|^{1/2}\bigr)^2\le 1 \le \bigl(|ad|^{1/2}+|bc|^{1/2}\bigr)^2.
$$
Taking square roots of both sides, including both choices of sign 
in the left hand inequality, gives (\ref{nn}), (\ref{nd}) and (\ref{cl1}).
\hfill$\square$

\medskip

\section{The proof of J{\o}rgensen's inequality}

\noindent
{\bf Proof of Theorem \ref{thmcp}. } 
Since (\ref{cond}) is invariant under conjugation, we may assume that $g$
is of the form (\ref{loxo}) and $h$ is  of the form (\ref{hform}).
Using (\ref{ratio1}) and (\ref{ratio2}) our hypothesis (\ref{cond}) 
can be rewritten as
\begin{equation}\label{rcond}
|ad|^{1/2}|bc|^{1/2}<\frac{1-M_g}{M_g^2}.
\end{equation}
Let $h_0 = h$  and  $h_{k+1} = h_kgh_k^{-1}$.   We write
$$
h_k=\left(
  \begin{array}{ccc}
    A_k& \alpha_k & \beta_k \\
   \eta_k & a_k& b_k\\
    \theta_k & c_k& d_k\\
    \end{array}
\right).
$$ Then
\begin{eqnarray*}
h_{k+1}&=& \left(
  \begin{array}{ccc}
     A_{k+1}& \alpha_{k+1} & \beta_{k+1} \\
    \eta_{k+1} & a_{k+1}& b_{k+1}\\
    \theta_{k+1} & c_{k+1}& d_{k+1}\\
    \end{array}
\right)\\
 &=& \left(
  \begin{array}{ccc}
    A_k& \alpha_k & \beta_k \\
   \eta_k & a_k& b_k\\
    \theta_k & c_k& d_k\\
    \end{array}
\right)\left(
       \begin{array}{ccc}
         L & 0 & 0\\
         0 & \lambda_n& 0 \\
         0& 0 & \overline{\lambda}_n^{-1} \\
             \end{array}
     \right)\left(
  \begin{array}{ccc}
     A_{k}^*& -\theta_k^*& -\eta_{k}^* \\
    -\beta_k^* & \overline{d_k}& \overline{b_k}\\
    -\alpha_k^* & \overline{c_k}& \overline{a_k}\\
    \end{array}
\right),
\end{eqnarray*}
where $L=diag(\lambda_1,\  \lambda_2,\  \cdots,\  \lambda_{n-1}).$
Therefore
 \begin{eqnarray}
\label{tolam1} 
a_{k+1} & = & -\eta_{k}L\theta_k^*+
a_k\lambda_n\overline{d_k}+b_k\overline{\lambda}_n^{-1}\overline{c_k},\\
\label{nn+1}
b_{k+1} & = & -\eta_kL\eta_k^*+
a_k\lambda_n\overline{b_k}+b_k\overline{\lambda}_n^{-1}\overline{a_k},\\
\label{n+1n}
c_{k+1} & = & -\theta_kL\theta_k^*+c_k\lambda_n\overline{d_k}
+d_k\overline{\lambda}_n^{-1}\overline{c_k},\\
\label{tolam2} 
d_{k+1} & = &  -\theta_kL\eta_k^*+c_k\lambda_n\overline{b_k}
+d_k\overline{\lambda}_n^{-1}\overline{a_k}.
\end{eqnarray}

\noindent
{\bf Claim 1:} 
We claim that if $|ad|^{1/2}|bc|^{1/2}<(1-M_g)/M_g^2$ then
$|b_kc_k|$ tends to $0$ as $k$ tends to infinity.

By (\ref{gg}) and (\ref{nn+1}), we have
 \begin{eqnarray*}
|b_{k+1}|&=&|\eta_k(I_{n-1}-L)\eta_k^*+
a_{k}(\lambda_n-1)\overline{b_{k}}+b_{k}(\overline{\lambda}_n^{-1}-1)\overline{a_{k}}|\\
&\leq
&\delta(g)\eta_k\eta_k^*+(|\lambda_n-1|+|\overline{\lambda}_n^{-1}-1|)|b_{k}a_{k}|\\
&=&
\delta(g)2\Re(a_{k}\overline{b_{k}})+(|\lambda_n-1|+|\overline{\lambda}_n^{-1}-1|)|b_{k}a_{k}|\\
&\leq
&(2\delta(g)+|\lambda_n-1|+|\overline{\lambda}_n^{-1}-1|)|b_{k}a_{k}|\\
&=& M_g|b_{k}a_{k}|.
\end{eqnarray*}

Similarly, by (\ref{tt})  and (\ref{n+1n})  we have
\begin{eqnarray*}
|c_{k+1}|=|\theta^{(k)}(I_{n-1}-L)\theta_k^*+
c_k(\lambda_n-1)\overline{d_k}+d_k(\overline{\lambda}_n^{-1}-1)\overline{c_k}|\leq
M_g|c_kd_k|.
\end{eqnarray*}

Therefore, for all $k\ge 0$ we have
\begin{equation}\label{eq-iteration}
|b_{k+1}c_{k+1}|^{1/2}\leq M_g|a_{k}d_k|^{1/2}|b_{k}c_k|^{1/2}.
\end{equation}
Using our hypothesis (\ref{rcond}) with $k=0$ this immediately gives
$$
|b_1c_1|^{1/2}\leq M_g|a_0d_0|^{1/2}|b_0c_0|^{1/2}<\frac{1-M_g}{M_g}.
$$
In particular,
$$
M_g\bigl(1+|b_1c_1|^{1/2}\bigr)<1.
$$
From this point on the proof closely follows the proof of the similar
result for complex hyperbolic space given by Jiang, Kamiya and Parker
\cite{jkp}.

We claim that for $k\ge 1$ we have
\begin{equation}\label{eq-ind-step}
|b_{k}c_{k}|^{1/2}\le 
\Bigl(M_g\bigl(1+|b_1c_1|^{1/2}\bigr)\Bigr)^{k-1}|b_1c_1|^{1/2}.
\end{equation}
In particular,
$$
|b_{k}c_{k}|^{1/2}\le |b_1c_1|^{1/2}.
$$
Certainly (\ref{eq-ind-step}) is true for $k=1$. Assume that 
(\ref{eq-ind-step}) is true for some $k\ge 1$.
Then, using (\ref{eq-iteration}) and (\ref{nn}), we have
\begin{eqnarray*}
|b_{k+1}c_{k+1}|^{1/2} & \le & M_g|a_{k}d_{k}|^{1/2}|b_{k}c_{k}|^{1/2} \\
& \le & M_g\bigl(1+|b_{k}c_{k}|^{1/2}\bigr)|b_{k}c_{k}|^{1/2} \\
& \le & M_g\bigl(1+|b_1c_1|^{1/2}\bigr)|b_{k}c_{k}|^{1/2} \\
& \le & M_g\bigl(1+|b_1c_1|^{1/2}\bigr)
\Bigl(M_g\bigl(1+|b_1c_1|^{1/2}\bigr)\Bigr)^{k-1}|b_1c_1|^{1/2} \\
& = & \Bigl(M_g\bigl(1+|b_1c_1|^{1/2}\bigr)\Bigr)^k|b_1c_1|^{1/2}.
\end{eqnarray*}
Then (\ref{eq-ind-step}) is true for $k+1$. The result follows by
induction.

Since $M_g\bigl(1+|b_1c_1|^{1/2}\bigr)<1$, an immediate consequence
of (\ref{eq-ind-step}) is that 
\begin{equation}\label{eq-limit-bc}
\lim_{k\rightarrow\infty}|b_kc_k|^{1/2}=0.
\end{equation}
This proves Claim 1.

\medskip

\noindent
{\bf Claim 2:} If  there exists some integer $k$ such that
\begin{equation}\label{dcond}
b_kc_k=0,
\end{equation}
then $\langle h,g\rangle$ is either elementary or not discrete.

If $b_k=0$ then, by (\ref{bb}), we have $\beta_k=0$ and $h_k(o)=o$.
Similarly, if $c_k=0$ then, by (\ref{aa}), we have $\alpha_k=0$
and so $h_k(\infty)=\infty$. If $b_kc_k=0$ but either $b_k$ or
$c_k$ is non-zero then $h_k$ fixes exactly one of $o$ and $\infty$.
Hence, $\langle g,\,h_k \rangle$ is not discrete by 
Theorem 3.1 of Kamiya \cite{kam83}. This implies that
$\langle g,\,h \rangle$ is not discrete.

Suppose then that $b_k=c_k=0$ for some $k\ge 1$. Then $h_k$ fixes both $o$ and 
$\infty$. In particular,
$$
o=h_k(o)=h_{k-1}gh_{k-1}^{-1}(o) \quad \hbox{ and }\quad
\infty=h_k(\infty)=h_{k-1}gh_{k-1}^{-1}(\infty).
$$
This means that $g$ fixes $h_{k-1}^{-1}(o)$ and $h_{k-1}^{-1}(\infty)$. 
If $k\ge 2$ then $h_{k-1}$ is loxodromic and so cannot swap $o$ and $\infty$.
Thus $h_{k-1}(o)=o$ and $h_{k-1}(\infty)=\infty$. By induction, we find that
$g$ fixes $h_0^{-1}(o)$ and $h_0^{-1}(\infty)$. In other words $h_0=h$ preserves
the set $\{o,\,\infty\}$ and so $\langle g,h\rangle$ is elementary.

This proves Claim 2.

\medskip

\noindent
{\bf Claim 3:}  If  
\begin{equation}\label{neq0}
\lim_{k\rightarrow\infty}|b_kc_k|=0\quad \hbox{ and }\quad
b_kc_k\neq 0 \quad \hbox{ for all } k\ge 1
 \end{equation}
then  $\langle h,g \rangle$ is  not discrete.

Assume that (\ref{neq0}) holds. Then from (\ref{nn}) we have
$$
|a_kd_k|^{1/2}\le |b_kc_k|^{1/2}+1
$$ 
and so $|a_kd_k|$ is bounded as $k$ tends to infinity.
Hence, from (\ref{ba}) and (\ref{gt}) we have
$$
|\beta_k^*\alpha_k| \leq 2|a_kd_k|^{1/2}|b_kc_k|^{1/2} \quad
\hbox{ and }\quad
|\eta_k\theta_k^*| \leq  2|a_kd_k|^{1/2}|b_kc_k|^{1/2}
$$
and so 
$$
\lim_{k\to\infty} |\beta_k^*\alpha_k| 
=\lim_{k\to\infty} |\eta_k\theta_k^*| =0.
$$
Likewise, 
$$
\lim_{k\to\infty} \eta_kL\theta_k^* =\lim_{k\to\infty}\theta_kL\eta_k^*=0.
$$
From (\ref{bba}) we have
$$
\lim_{k\to\infty} \overline{d}_ka_k=
\lim_{k\to\infty} \bigl(1+\beta_k^*\alpha_k-\overline{b}_kc_k\bigr)=1.
$$
Therefore, from (\ref{tolam1}) and (\ref{tolam2}) we have
\begin{eqnarray}
\label{lima}\lim_{k\to\infty} |a_{k+1}| 
& = & \lim_{k\to\infty}
\bigl|-\eta_{k}L\theta_k^*+
a_k\lambda_n\overline{d_k}+b_k\overline{\lambda}_n^{-1}\overline{c_k}\bigr|
\ =\ |\lambda_n|, \\
\label{limd}\lim_{k\to\infty} |d_{k+1}|
& = & \lim_{k\to\infty} \bigl| -\theta_kL\eta_k^*+c_k\lambda_n\overline{b_k}
+d_k\overline{\lambda}_n^{-1}\overline{a_k}\bigr|
\ =\ |\lambda_n|^{-1}.
\end{eqnarray}

When proving Claim 1, we showed that
$$
|b_{k+1}|\le M_g|a_k|\,|b_k|\quad \hbox{ and } \quad
|c_{k+1}|\le M_g|d_k|\,|c_k|.
$$
Since $M_g<1$ we can find $K$ so that, using (\ref{lima}) and
(\ref{limd}), for all $k\ge K$
$$
M_g|a_k|<|\lambda_n|\quad\hbox{ and } \quad M_g|d_k|<|\lambda_n|^{-1}.
$$
Hence there exist constants $\kappa_1$ and $\kappa_2$ so that, for all $k\ge K$
$$
M_g|a_k|\,|\lambda_n|^{-1}<\kappa_1<1 \quad\hbox{ and } \quad 
M_g|d_k|\,|\lambda_n|<\kappa_2<1.
$$
Therefore for $k\ge K$
\begin{eqnarray*}
|b_{k+1}|\,|\lambda_n|^{-k-1}
& \le &  \bigl(M_g|a_k|\,|\lambda_n|^{-1}\bigr)|b_k|\,|\lambda_n|^{-k}
\ < \ \kappa_1|b_k|\,|\lambda_n|^{-k}
\ \le \  \kappa_1^{k+1-K}|b_K|\,|\lambda_n|^{-K}, \\
|c_{k+1}|\,|\lambda_n|^{k+1}
& \le &  \bigl(M_g|d_K|\,|\lambda_n|\bigr)|c_k|\,|\lambda_n|^{k}
\ <\ \kappa_2|c_k|\,|\lambda_n|^{k}
\ \le \  \kappa_2^{k+1-K}|c_K|\,|\lambda_n|^{K}.
\end{eqnarray*}
Since $K$ was chosen so that $\kappa_i<1$ for $i=1,\,2$, we see that
\begin{equation}\label{limbc}
\lim_{k\to\infty}|b_k|\,|\lambda_n|^{-k}=0\quad \hbox{ and }\quad
\lim_{k\to\infty}|c_k|\,|\lambda_n|^{k}=0.
\end{equation}

Following J{\o}rgensen, we now define the sequence
$f_k=g^{-k}h_{2k}g^{k}$.   As a matrix in $Sp(n,1)$ this is given by
\begin{equation}f_k=\left(\begin{matrix}
L^{-k}A_{2k}L^k & L^{-k}\alpha_{2k}\lambda_n^k & 
L^{-k}\beta_{2k}\overline{\lambda}_n^{-k} \\
\lambda_n^{-k}\eta_{2k}L^k & \lambda_n^{-k}a_{2k}\lambda_n^{k} & 
\lambda_n^{-k}b_{2k}\overline{\lambda}_n^{-k} \\
\overline{\lambda}_n^k\theta_{2k}L^k& 
\overline{\lambda}_n^kc_{2k}\lambda_n^{k} & 
\overline{\lambda}_n^kd_{2k}\overline{\lambda}_n^{-k} 
\end{matrix}\right).
\end{equation}
Using (\ref{lima}) and (\ref{limd}), we have
$$
\lim_{k\to\infty}|\lambda_n^{-k}a_{2k}\lambda_n^{k}|
=\lim_{k\to\infty}|a_{2k}|=|\lambda_n| \quad \hbox{ and }\quad
\lim_{k\to\infty}|\overline{\lambda}_n^kd_{2k}\overline{\lambda}_n^{-k}|
=\lim_{k\to\infty}|d_{2k}|=|\lambda_n|^{-1}.
$$
Similarly, using (\ref{limbc}). we have
$$
\lim_{k\to \infty}|\lambda_n^{-k}b_{2k}\overline{\lambda}_n^{-k}|
=\lim_{k\to\infty}|b_{2k}|\,|\lambda|^{-2k}=0
\quad \hbox{ and }\quad
\lim_{k\to\infty}|\overline{\lambda}_n^kc_{2k}\lambda_n^{k}|
=\lim_{k\to\infty}|c_{2k}|\,|\lambda|^{2k}=0.
$$
Then, using (\ref{aa}), (\ref{bb}), (\ref{gg}) and (\ref{tt}) for the matrix 
$f_k$, we have
\begin{eqnarray*}
\lim_{k\to\infty}|L^{-k}\alpha_{2k}\lambda_n^k|^2
& \le & \lim_{k\to\infty}2|\overline{\lambda}_n^kc_{2k}\lambda_n^{k}|\,
|\lambda_n^{-k}a_{2k}\lambda_n^{k}| \ = \ 0, \\
\lim_{k\to\infty}|L^{-k}\beta_{2k}\overline{\lambda}_n^{-k}|^2
& \le & \lim_{k\to\infty}2|\lambda_n^{-k}b_{2k}\overline{\lambda}_n^{-k}|\,
|\overline{\lambda}_n^kd_{2k}\overline{\lambda}_n^{-k}|\ = \ 0, \\
\lim_{k\to\infty}|\lambda_n^{-k}\eta_{2k}L^k |^2
& \le & \lim_{k\to\infty}2|\lambda_n^{-k}b_{2k}\overline{\lambda}_n^{-k}|\,
|\lambda_n^{-k}a_{2k}\lambda_n^{k}| \ = \ 0, \\
\lim_{k\to\infty}|\overline{\lambda}_n^k\theta_{2k}L^k|^2
& \le & \lim_{k\to\infty}2|\overline{\lambda}_n^kc_{2k}\lambda_n^{k}|\,
|\overline{\lambda}_n^kd_{2k}\overline{\lambda}_n^{-k}|\ = \ 0.
\end{eqnarray*}
Finally, this means that $L^{-k}\alpha_{2k}\lambda_n^k$ and 
$L^{-k}\beta_{2k}\overline{\lambda}_n^{-k}$ both tend to the zero vector.
Hence, using (\ref{AA}) on the matrix $f_k$, we see that
\begin{eqnarray*}
\lefteqn{\lim_{k\to\infty} (L^{-k}A_{2k}L^k)(L^{-k}A_{2k}L^k)^*} \\
& = & I_{n-1}+\lim_{k\to\infty}\Bigl(
(L^{-k}\alpha_{2k}\lambda_n^k)
(L^{-k}\beta_{2k}\overline{\lambda}_n^{-k})^*
+(L^{-k}\beta_{2k}\overline{\lambda}_n^{-k})
(L^{-k}\alpha_{2k}\lambda_n^k)^*\Bigr) \\
& = & I_{n-1}.
\end{eqnarray*}
Therefore $\{f_k\ :\ k\ge K\}$ lies in a compact subset of ${\rm Sp}(n,1)$
and so contains a convergent subsequence.
This proves Claim 3, and hence completes the proof of Theorem \ref{thmcp}.

\hfill$\square$

\medskip

\noindent
{\bf Proof of Corollary \ref{corcp}. }
Without loss of
generality, we assume $u=\infty$ and $v=o$, and $g$ is of  the form
(\ref{loxo}) and $h$ is  of the form (\ref{hform}).
Using the identities (\ref{ratio1}), (\ref{ratio2}) and (\ref{ratio3}) 
from Lemma \ref{lem-X-ratio}, the conditions
(\ref{cond2a}), (\ref{cond2b}), (\ref{cond-1.1}) and 
(\ref{condK2-improved}) can be  rewritten as
\begin{eqnarray}
|bc|^{1/2} & < & \frac{1-M_g}{M_g}, \label{cond2a-1} \\
|ad|^{1/2} & < & \frac{1-M_g}{M_g}, \label{cond2b-1} \\
\frac{|bc|^{1/2}}{|ad|^{1/2}} & < &  1-M_g, \label{cond-1.1-1} \\
|ad|+|bc| & < & \frac{2(1-M_g)}{M_g^2}. \label{condK2-improved-1}
\end{eqnarray}
Our strategy will be to show that each of these conditions
implies (\ref{rcond}) and the result will then follow from Theorem \ref{thmcp}.

Using (\ref{nn}) condition (\ref{cond2a-1}) implies
$$
|ad|^{1/2}|bc|^{1/2}\le\bigl(1+|bc|^{1/2}\bigr)|bc|^{1/2}
< \left(1+\frac{1-M_g}{M_g}\right)\frac{1-M_g}{M_g}=\frac{1-M_g}{M_g^2}.
$$
Similarly, using (\ref{nd}), condition (\ref{cond2b-1}) gives
(\ref{rcond}). 

Using (\ref{nn}) condition (\ref{cond-1.1-1}) implies
$$
M_g < 1-\left|\frac{bc}{ad}\right|^{1/2} 
\le 1-\frac{|bc|^{1/2}}{|bc|^{1/2}+1}
= \frac{1}{|bc|^{1/2}+1}.
$$
Rearranging, this is equivalent to (\ref{cond2a-1}) and so the result 
follows from the earlier part of this proof.

Finally, condition (\ref{condK2-improved-1}) implies
$$
|ad|^{1/2}|bc|^{1/2}\le \frac{1}{2}\bigl(|ad|+|bc|\bigr)
\le \frac{1-M_g}{M_g^2}.
$$
Therefore in each case $\langle h,g\rangle$ is either elementary or 
not discrete by Theorem \ref{thmcp}.
\hfill$\square$

\medskip

\section{Collars in ${\bf H}_{\bh}^n$}

We need the following lemma, whose proof can be  verified
directly, to prove Theorem \ref{nthm1.4}.

\begin{lem}
Let ${\bf p},\,{\bf q}\in V_0$ be null vectors with 
$\langle {\bf p}, {\bf q} \rangle=-1.$ For all real $t$ let 
$\gamma(t)$ be the point in ${\bf H}_{\bh}^n$ corresponding to the vector 
$e^{\frac{t}{2}}{\bf p}+e^{-\frac{t}{2}}{\bf q}$ in $\bh^{n,1}$. Then
$\gamma=\{\gamma(t)|t\in \br\}$ is the geodesic in ${\bf H}_{\bh}^n$ with
endpoints $\bp({\bf p})$ and  $\bp({\bf q})$ parametrised by arc length
$t$.
\end{lem}

\medskip

The following Proposition relates cross-ratios to 
the distance between geodesics. It will be crucial in our proofs of
Theorems \ref{nthm1.4} and \ref{nthm1.6}

\begin{prop}\label{prop-dist}
Let $\gamma_1$ and $\gamma_2$ be geodesics in ${\bf H}_{\bh}^n$ with
endpoints $u_1$, $v_1$ and $u_2$, $v_2$ respectively. Then
$$
\cosh\bigl(\rho(\gamma_1,\gamma_2)\bigr)
\geq \bigl|[v_2,u_1,v_1,u_2]\bigr|+\bigl|[v_2,v_1,u_1,u_2]\bigr|. 
$$
\end{prop}

\noindent
{\bf Proof.} 
Without loss of generality, suppose that $u_1=o$ and $v_1=\infty$.
Also, let $h\in {\rm PSp}(n,1)$ be a map so that $u_2=h(o)$ and
$v_2=h(\infty)$. Suppose that $h\in G$ has the form (\ref{hform}), and so the 
cross-ratios are given by (\ref{ratio1}) and (\ref{ratio2}).
Let $p_t$ and $q_s$ be two points on the geodesic $\gamma_1$ 
and $\gamma_2=h(\gamma_1)$, respectively.  Then, letting ${\bf 0}$ denote 
the zero vector in $\bh^{n-1}$, we can choose $t,\,s\in \br$ such that
$$
{\bf p}_t=
\left(\begin{matrix} {\bf 0} \\ e^{-t} \\ 1 \end{matrix}\right)
\in \bp^{-1}(\gamma) ,\quad
{\bf q}_s=h({\bf p}_s)
=\left(\begin{matrix} A & \alpha & \beta \\ \eta & a & b \\ \theta & c & d
\end{matrix}\right)
\left(\begin{matrix} {\bf 0} \\ e^{-s} \\ 1 \end{matrix}\right)
=\left(\begin{matrix} \alpha e^{-s}+\beta \\ 
ae^{-s}+b \\ ce^{-s}+d \end{matrix}\right) \in \bp^{-1}(\gamma).
$$
Since $\langle {\bf p}_t,\,{\bf p}_t\rangle=-2e^{-t}$,  
$\langle{\bf q}_s,\,{\bf q}_s\rangle= \langle h({\bf p}_s),\,h({\bf
p}_s)\rangle=\langle {\bf p}_s,\,{\bf p}_s\rangle=-2e^{-s}$ and 
$$
\langle {\bf p}_t,\,{\bf q}_s\rangle \langle {\bf q}_s,\,{\bf p}_t\rangle
=\Bigl(\overline{a}e^{-s}+\overline{b}+(\overline{c}e^{-s}
+\overline{d})e^{-t}\Bigr)\Bigl(ae^{-s}+b+(ce^{-s}+d)e^{-t}\Bigr),
$$  
we have
\begin{eqnarray*}
\cosh\bigl(\rho(p_t,q_s)\bigr) & = & 2\cosh^2\frac{\rho(p_t,q_s)}{2} -1 \\
&=& 2\,\frac{\langle {\bf p}_t,\,{\bf q}_s\rangle
\langle {\bf q}_s,\,{\bf p}_t\rangle}
{\langle {\bf p}_t,\,{\bf p}_t\rangle
\langle{\bf q}_s,\,{\bf q}_s\rangle} -1 \\
&=&\frac{1}{2}\Bigl(|a|^2e^{t-s}+|d|^2e^{s-t}+|c|^2e^{-(s+t)}+|b|^2e^{s+t}
+(\overline{d}a+\overline{b}c)+(\overline{a}d+\overline{c}b)\\
&&\quad +(\overline{a}b+\overline{b}a)e^{t}+(\overline{c}d+\overline{d}c)e^{-t}
+(\overline{d}b+\overline{b}d)e^{s}+(\overline{c}a+\overline{a}c)e^{-s}-2\Bigr).
\end{eqnarray*}
By (\ref{AA})-(\ref{aa}) and the property
$a\overline{b}+b\overline{a}=\overline{b}a+\overline{a}b$ for $a,b \in \bh$,  we have
\begin{eqnarray*}
\cosh\bigl(\rho(p_t,q_s)\bigr)
&=&\frac{1}{2}\Bigl(|a|^2e^{t-s}+|d|^2e^{s-t}+|c|^2e^{-(s+t)}+|b|^2e^{s+t}
+\beta^*\alpha+\alpha^*\beta\\
&&\quad +\eta\eta^*e^{t}+\theta\theta^*e^{-t}
+\beta^*\beta e^{s}+\alpha^*\alpha e^{-s}\Bigr) \\
&\geq& |ad|+|bc|+\Re(\beta^*\alpha)+|\eta|\,|\theta|
+|\beta^*\alpha| \\
&\geq& |ad|+|bc| \\
&=&\bigl|[v_2,u_1,v_1,u_2]\bigr|+\bigl|[v_2,v_1,u_1,u_2]\bigr|. 
\end{eqnarray*}
This is true for all points $p_t$ and $q_s$ and so it 
proves the proposition.
\hfill$\square$

\medskip

\noindent
{\bf Proof of Theorem \ref{nthm1.4}.}\quad  
Without loss of
generality, we suppose that $g$ has the form (\ref{loxo}) and so fixes 
$o$ and $\infty$. If $h\in G$ maps $\gamma$ to itself then it must
map $T_r(\gamma)$ to itself. 

Therefore we suppose that $h$ does not map $\gamma$ to itself. We must
show that $T_r(\gamma)$ is disjoint from its image under $h$.
We first use Proposition \ref{prop-dist} to estimate the distance between
$\gamma$ and $h(\gamma)$ and then use condition (\ref{condK2-improved})
from Corollary \ref{corcp} to conclude that, since $G$ is discrete
and non-elementary, we have
\begin{eqnarray*}
\cosh\bigl(\gamma,h(\gamma)\bigr)
& \ge & \bigl|[h(\infty),o,\infty,h(o)]\bigr| 
+\bigl|[h(\infty),\infty,o,h(o)]\bigr| \\
& \ge & \frac{2(1-M_g)}{M_g^2}.
\end{eqnarray*}
This implies that $T_r(\gamma)$ is disjoint from its image under $h$.
\hfill$\square$

\medskip

\noindent
{\bf Proof of Theorem \ref{nthm1.6}.} Let
$\mathcal{M}={\bf H}_{\bh}^n/G$ where $G$ is a discrete, non-elementary,
torsion-free subgroup of ${\rm Sp}(n,1)$.  Let $g$ and $h$ be two
loxodromic elements of $G$ whose axes, $\gamma_1$ and $\gamma_2$,
project to distinct short, simple, closed geodesics
$\gamma_1'=\gamma_1/\langle g \rangle$.   and
$\gamma_2'=\gamma_2/\langle h \rangle$.  Reordering if necessary,
suppose that $M_h\le M_g$. Consider tubes $T_{r_1}(\gamma_1)$ and
$T_{r_2}(\gamma_2)$ around $\gamma_1$ and $\gamma_2$ where
$$
\cosh(2r_1)=\frac{2(1-M_g)}{M_g^2}
 \quad \hbox{ and } \quad
\cosh(2r_2)=\frac{2(1-M_h)}{M_h^2}.
$$
We want to show that these tubes are disjoint. It suffices to show that
$\rho(\gamma_1,\gamma_2)\ge r_1+r_2$.

Without loss of generality,
we suppose that  $g$ is of the form  (\ref{loxo}) and  $h$ has fixed
points $p=(p_1,\cdots, p_n)^t\in \partial {\bf H}_{\bh}^n$ and
$q=(q_1,\cdots, q_n)^t\in \partial {\bf H}_{\bh}^n.$  That is
$$
\sum_{i=1}^{n-1}{|p_i|^2}=p_n+\overline{p}_n, \ \
\sum_{i=1}^{n-1}{|q_i|^2}=q_n+\overline{q}_n.
$$
Let ${\bf p}=(p_1,\cdots, p_n,1)^t\in \bp^{-1}(p)$ and 
${\bf q}=(q_1,\cdots, q_n,1)^t\in \bp^{-1}(q)$. 
Then by the definition of quaternionic cross-ratio, we have
$$
\bigl|[o,q,p,\infty]\bigr|
=\frac{|p_n|}{|\langle {\bf p},{\bf q}\rangle |},\qquad 
\bigl|[o,p,q,\infty]\bigr|
=\frac{|q_n|}{|\langle {\bf p},{\bf q} \rangle|}.
$$
Direct computation implies that
\begin{eqnarray*}
\bigl|\langle g({\bf q}),{\bf q}\rangle\bigr|
&=&\left|-(\overline{q}_n\overline{\lambda}_n^{-1}+\lambda_n q_n)
+\sum_{i=1}^{n-1}{\overline{q}_i\lambda_i q_i}\right|\\
&=&\left|(1-\lambda_n)q_n+\overline{q}_n(1-\overline{\lambda}_n^{-1})
+\sum_{i=1}^{n-1}{\overline{q}_i(\lambda_i-1) q_i}\right|\\
&\leq&|\lambda_n-1||q_n|+|\overline{\lambda}_n^{-1}-1||q_n|
+\delta(g) \sum_{i=1}^{n-1}|q_i|^2 \\
&=&|\lambda_n-1||q_n|+|\overline{\lambda}_n^{-1}-1||q_n|
+\delta(g)(q_n+\overline{q}_n) \\
&\leq & M_g|q_n|.
\end{eqnarray*}
Similarly, we have $|\langle g({\bf p}),{\bf p}\rangle|\leq  M_g|p_n|$.
We remark that these equations are special cases of (\ref{eq-mp}). 
Therefore, we get
\begin{equation}
M_g^2\bigl|[o,p,q,\infty]\bigr|\,\bigl|[o,q,p,\infty]\bigr|
=\frac{M_g|q_n|}{|\langle{\bf p},{\bf q}\rangle|}\,
\frac{M_g|p_n|}{|\langle{\bf p},{\bf q}\rangle|}
\geq
\frac{|\langle g({\bf q}),{\bf q}\rangle||\langle g({\bf p}),{\bf p}\rangle|}
{|\langle{\bf p},{\bf q}\rangle|^2}
=\bigl|[g(p),q,p,g(q)]\bigr| .
\end{equation}
Using Proposition \ref{prop-dist}, we get
\begin{eqnarray*}
\cosh\bigl(\rho(\gamma_1,\gamma_2)\bigr)
&\geq&\bigl|[o,p,q,\infty]\bigr|+\bigl|[o,q,p,\infty]\bigr| \\
&\geq& 2\bigl|[o,p,q,\infty]\bigr|^{1/2}\bigl|[o,q,p,\infty]\bigr|^{1/2} \\
&\geq& \frac{2}{M_g}\bigl|[g(p),q,p,g(q)]\bigr|^{1/2}.
\end{eqnarray*}
Using condition (\ref{cond2a}) with the roles of $g$ and $h$ interchanged 
we have $\bigl|[g(p),q,p,g(q)]\bigr|^{1/2}\ge (1-M_h)/M_h$.
Using this and $M_g\ge M_h$ we have
\begin{eqnarray*}
\cosh^2\bigl(\rho(\gamma_1,\gamma_2)\bigr)
& \geq & \left(\frac{2(1-M_h)}{M_gM_h}\right)^2 \\
& \geq & \left(\frac{2(1-M_g)}{M_g^2}\right)\left(\frac{2(1-M_h)}{M_h^2}\right)
 \\
&=&\cosh(2r_1)\cosh(2r_2)\\
&\geq&\cosh^2(r_1+r_2).
\end{eqnarray*}
Therefore $\rho(\gamma_1,\gamma_2)\ge r_1+r_2$ as required.
\hfill$\square$

\medskip

\section{Collar width solely in terms of geodesic length}

In the complex case Markham and Parker \cite{mapa}
used a lemma of Zagier to obtain the width of the tubular neighbourhood of
the simple geodesic $\gamma$ entirely in terms of its length.
In this section, we will consider the counterpart in $n$-dimensional
quaternionic hyperbolic manifold. First, we give an example to
illustrate our idea.

\begin{example}\label{example1} 
Let $G$ be a discrete, torsion-free, non-elementary subgroup of
${\rm Sp}(2,1)$ with
$$
g={\rm diag}(e^{{\bf i}\beta},e^{l/2+{\bf i}\alpha},e^{-l/2+{\bf i}\alpha})\in G.
$$ 
(Here the imaginary units that generate ${\bh}$ are denoted ${\bf i}$, 
${\bf j}$ and ${\bf k}$ in order to distinguish them from the
indices denoted by $i$, $j$ and $k$.)
Define $f(k)=M_{g^k}=\bigl|e^{kl/2+{\bf i}k\alpha}-1\bigr|
+\bigl|e^{-kl/2+{\bf i}k\alpha}-1\bigr|+2\bigl|e^{{\bf i}k\beta}-1\bigr|$. 
Then we have
\begin{equation}\label{poi-2}
f(k)=2\sqrt{\left(\cosh\frac{kl}{2}+1\right)\left(\cosh
\frac{kl}{2}-\cos(k\alpha)\right)}+2\sqrt{2\Bigl(1-\cos(k\beta)\Bigr)}.
\end{equation}
Consider the case where
$l=10^{-3}, \alpha=\frac{\pi}{3}, \beta=\frac{\pi}{4}$.

We see, for $k\in {\mathbb Z}$, that if $k$ is not a multiple of $8$ then
$\cos(k\beta)\le 1/\sqrt{2}$ and so
$$
f(k)\ge 2\sqrt{2\Bigl(1-\cos(k\pi/4)\Bigr)}
\ge 2\sqrt{2-\sqrt{2}}>1.
$$
Likewise, when $k$ is not a multiple of $6$ then $\cos(k\alpha) \le 1/2$
and so
$$
f(k)\ge 2\sqrt{\left(\cosh\frac{k}{2000}+1\right)
\left(\cosh\frac{k}{2000}-\cos\frac{k\pi}{3}\right)}
\ge 2\sqrt{2\left(1-\cos\frac{k\pi}{3}\right)}\ge 2\sqrt{2-1}=2.
$$
On the other hand, if $k$ is a multiple of both $8$ and $6$, that is
a multiple of $24$, then
$$
f(k) = 2\sqrt{\left(\cosh\frac{k}{2000}+1\right)\left(\cosh
\frac{k}{2000}-1\right)}
=2\sinh\frac{k}{2000}.
$$
Hence as $k$ ranges over positive integers, the minimum value of
$f(k)$ is attained for $k=24$ and is approximately 
$0.024$.
\end{example}

\medskip

The above example shows that when $l$ and $k$ are small, then
$\cos(k\alpha)$ and $\cos(k\beta)$ contribute  the dominant part of
the value $f(k)$.  Although $f(k)\to \infty$ as $k\to \infty$, we
sometimes can choose suitable $k$ such that $k\alpha$ and $k\beta$
are close to multiples of $2\pi$ which may lead to
$f(k)<\sqrt{3}-1$. This observation gives an improvement of Theorem
\ref{nthm1.4} by replacing $M_g$ with a suitable $M_{g^k}$.

We now investigate how $M_{g^k}$ varies with $k$.
Let $g$ be of the form (\ref{loxo}). We can conjugate all its right
eigenvalues to unique complex numbers with non-negative imaginary
part, that is,
\begin{eqnarray*}
\lambda_i & = & u_ie^{{\bf i}\beta_i}u_i^{-1}, \hbox{ for } 1 \leq i\leq n-1,\\
\lambda_n & = & u_n|\lambda_n|e^{{\bf i}\beta_n}u_n^{-1},\\
\overline{\lambda}_n^{-1} & = & u_{n+1}|\lambda_n|^{-1}e^{{\bf i}\beta_n}u_{n+1}^{-1},
\end{eqnarray*}
where $0\leq \beta_i \leq \pi$ for $1\le i\le n$.
Recall that
\begin{eqnarray*}
M_g & = & \bigl|\lambda_n-1\bigr|+\bigl|\overline{\lambda}_n^{-1}-1\bigr|
+2\max_{1\le i\le n-1}|\lambda_i-1| \\
& = & 2\sqrt{\left(\cosh\frac{l}{2}+1\right)
\left(\cosh\frac{l}{2}-\cos(\beta_n)\right)}
+\max_{1\le i\le n-1}2\sqrt{2\Bigl(1-\cos(\beta_i)\Bigr)}.
\end{eqnarray*}
Since the eigenvalues of $M^k$ are $\lambda_i^k$, $\lambda_n^k$, 
$\overline{\lambda}_n^{-k}$ we have
\begin{equation}\label{poi}
M_{g^k} = 2\sqrt{\left(\cosh\frac{kl}{2}+1\right)\left(\cosh
\frac{kl}{2}-\cos(k\beta_n)\right)}
+\max_{1\le i\le n-1}2\sqrt{2\Bigl(1-\cos(k\beta_i)\Bigr)}.
\end{equation}
Define $T$ to be the minimum value of $M_{g^k}$. That is
\begin{equation}\label{min}
T= \min_{1\leq k <\infty}\left\{
2\sqrt{\left(\cosh\frac{kl}{2}+1\right)\left(\cosh
\frac{kl}{2}-\cos(k\beta_n)\right)}
+\max_{1\le i\le n-1}2\sqrt{2\Bigl(1-\cos(k\beta_i)\Bigr)}\right\}.
 \end{equation}
Then by Theorem \ref{nthm1.4} we have the following corollary.

\begin{cor} \label{coo2}Let $G$ be a discrete, non-elementary, torsion-free
subgroup of ${\rm Sp}(n,1)$.  Let $g$ be a loxodromic element of $G$ with
axis the geodesic $\gamma$.  Let $T$ be given by (\ref{min}) and suppose
that $T <\sqrt{3}-1$. Let $r$ be positive real number defined by
$$
\cosh(2r)=\frac{2(1-T)}{T^2}.
$$
Then the tube $T_r(\gamma)$ is precisely invariant under $G$.
\end{cor}

In order to prove Theorem \ref{nthm1.7}, we need the following
adapted   Pigeonhole Principle.

\begin{lem} \label{meyer} (cf. Pigeonhole Principle in \cite{mey})
Given  $0\leq \beta_1,\cdots,\beta_n <2\pi$ and a positive integer
$N\geq 2$, there exists $k\leq N^n$ such that
$$
\cos(k\beta_i)\geq \cos\frac{2\pi}{N},\quad
$$
for each $1\leq i\leq n$.
\end{lem}

\noindent
{\bf Proof.} Consider the solid $n$-cube $[0,2\pi]^n$ in $\br^n$ and
consider the points $z_k=(k\beta_1,k\beta_2...,k\beta_n)$ for
each $1\le k\le N^n.$  There are $N^n$ of them. For each
$i\in\{1,\,\ldots,\,n\}$ and each $k\in\{1,\,\ldots,\,N^n\}$ let
$m_{ik}$ be an integer so that $k\beta_i-2\pi m_{ik}\in[0,2\pi)$.
Let
$$
\hat z_k=\bigl(k\beta_1-2\pi m_{1k},\,\ldots,\,k\beta_n-2\pi m_{nk}\bigr)
\in[0,2\pi]^n.
$$
Divide the $n$-cube into $N^n$ cubes of side length $2\pi/N$ and consider 
which of these small cubes contain the points $\hat z_k$.

If, for some $j\in\{1,\,\ldots,\,N^n\}$, the point $\hat z_j$ 
lies in the $n$-cube $[0,\frac{2\pi}{N}]^n$ then
$\cos(j\beta_i)\ge \cos\frac{2\pi}{N}$ for each $i\in\{1,\,\ldots,\, n\}$
and we have the result.

Suppose that none of the points $\hat z_k$ lie in the $n$-cube
$[0,\frac{2\pi}{N}]^n$. Then there is at least one small $n$-cube
with two points in it, say $\hat z_j$ and $\hat z_k$,  where $j>k.$
Then $\hat z_{j-k}$ is in the $n$-cube $I_1\times I_2\times
\cdots\times I_n$, where each $I_i=[0,\frac{2\pi}{N}]$ or
$[2\pi-\frac{2\pi}{N},2\pi].$   That is
$
\cos\bigl((j-k)\beta_i\bigr)\ge \cos\frac{2\pi}{N}
$
for each $i\in\{1,\,\ldots,\, n\}$.
The proof is complete.
\hfill$\square$

\medskip

\begin{figure}
\begin{center}
  \includegraphics[width=12cm, height=6cm]{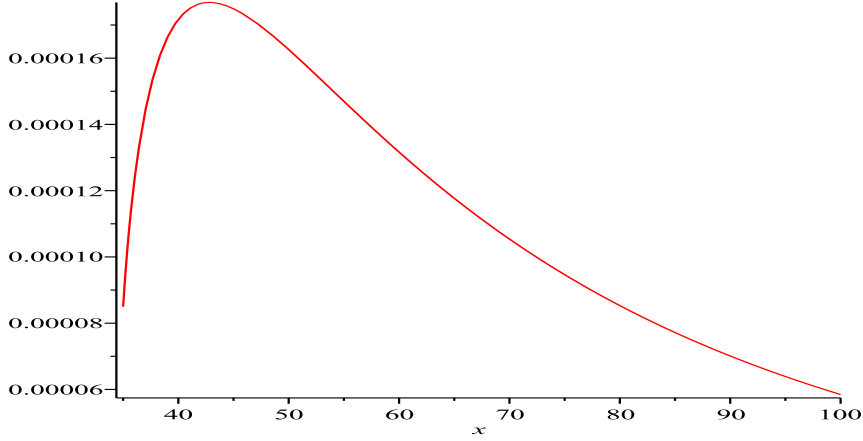}\\
  \caption{
The graph of function $l=l(x)$ defined by (\ref{hhl1})
for $n=2$.}\label{fig1}
  \end{center}
\end{figure}

\noindent
{\bf Proof of Theorem \ref{nthm1.7}.} As in (\ref{poi}), we have
$$
M_{g^k} = 2\sqrt{\left(\cosh\frac{kl}{2}+1\right)\left(\cosh
\frac{kl}{2}-\cos(k\beta_n)\right)}
+\max_{1\le i\le n-1}2\sqrt{2\Bigl(1-\cos(k\beta_i)\Bigr)}.
$$
By Lemma \ref{meyer} for $N\geq 2$, there exists $k\leq N^n$ such that
$$
\cos(k\beta_i)\ge \cos\frac{2\pi}{N}
$$
for each $1\le i\le n$.
Then there exists $k\leq N^n$ such that
\begin{equation}\label{fk-bound}
M_{g^k}\le 2\sqrt{\left(\cosh\frac{N^nl}{2}+1\right)
\left(\cosh\frac{N^nl}{2}-\cos\frac{2\pi}{N}\right)}+
2\sqrt{2\left(1-\cos\frac{2\pi}{N}\right)}.
\end{equation}
Define
$$
h(x,l)=2\sqrt{\left(\cosh\frac{x^nl}{2}+1\right)
\left(\cosh\frac{x^nl}{2}-\cos\frac{2\pi}{x}\right)}+
2\sqrt{2\left(1-\cos\frac{2\pi}{x}\right)}.
$$
Notice that $h(x,l)$ is an increasing function of $l$ for fixed $x$.
Therefore when $l>0$ we have
$$
h(x,l)> h(x,0)=4\sqrt{2\left(1-\cos\frac{2\pi}{x}\right)}.
$$
When $x\ge 2$ the function $h(x,0)$ is a decreasing function of $x$. 
Define $x_0$ by $h(x_0,0)=\sqrt{3}-1$, that is
$$
x_0= \frac{2\pi}{\arccos\frac{14+\sqrt{3}}{16}}\approx 34.284.
$$
Then if $2\le x < x_0$ and $l>0$ we have 
$h(x,l)>h(x,o)>h(x_0,0)=\sqrt{3}-1$.
Hence, in order to have $h(x,l) \le \sqrt{3}-1$ we must have $x \ge x_0$.
For all $x\ge x_0$ the equation $h(x,l)=\sqrt{3}-1$ defines a function
\begin{equation}\label{hhl1}
l(x)=\frac{2}{x^n}{\rm arccosh}\left(
\frac{\sqrt{13-2\sqrt{3}-6\cos\frac{2\pi}{x}+\cos^2\frac{2\pi}{x}
-4\bigl(\sqrt{3}-1\bigr)
\sqrt{2\left(1-\cos\frac{2\pi}{x}\right)}}-1+\cos\frac{2\pi}{x}}{2}\right).
\end{equation}
Hence for all integers $N\geq 35$, we can find $l$ satisfying the
condition (\ref{condi}). Then our result follows from the
application of Theorem  \ref{nthm1.4}. The proof is complete.
\hfill$\square$

\medskip

With the aid of mathematical software, for case $n=2$,  we find that
when $N=43$, we get the maximal interval $0<l<l(43)\approx
0.00017681$ to apply our theorem. The graph of function $l(x)$
defined by (\ref{hhl1}) is given in Figure \ref{fig1}. This gives
the proof of Corollary \ref{cor-1.8}.

Given the rotational angles of loxodromic element, we may be able to
use Corollary \ref{coo2} to choose suitable $N$ which may less than
35.  For instance in Example \ref{example1}, the optimum value
occurs when $N=24$.


\begin{thebibliography}{99}

\bibitem{bm} Basmajian, A., Miner, R.:  Discrete subgroups of complex
hyperbolic motions.   Invent. Math. {\bf 131}, 85-136 (1998)


\bibitem{caoh} Cao, W. S.:  J{\o}rgensen's inequality for 
quaternionic hyperbolic space. 
preprint(2008), available at  arXiv: 0906.1475v1

\bibitem{cpw}Cao, W. S., Parker, J. R., Wang, X. T.:  On the classification of
quaternionic M\"obius transformations.   Math. Proc. Cambridge
Philos. Soc. {\bf 137}, 349-361 (2004)





 \bibitem{caog} Cao, W. S., Tan, H . O.:  J{\o}rgensen's
inequality for quaternionic hyperbolic space with elliptic elements,
  Bull. Austral. Math. Soc. (to appear).






\bibitem{chen} Chen, S. S., Greenberg, L.:  Hyperbolic spaces, Contributions
to analysis.   Academic Press, New York.  49-87, (1974)

\bibitem{gold} Goldman, W. M.: Complex Hyperbolic Geometry.  Oxford  University Press
(1999)

 \bibitem{jkp} Jiang, Y. P., Kamiya, S., Parker, J. R.:
 J{\o}rgensen's inequality for complex hyperbolic space.   Geometriae Dedicata. {\bf
97}, 55-80  (2003)


 \bibitem{jor} J{\o}rgensen, T.:  On discrete groups of M\"obius transformations.  Amer. J. Math. {\bf
 98}, 739-749 (1976)

\bibitem{kam83} S. Kamiya.  Notes on non-discrete subgroups of U(1,n; F).
{\sl Hiroshima Math. J.} {\bf 13}  (1983), 501-506.



\bibitem{kim} Kim, D.:  Discreteness criterions of isometric subgroups for quaternionic
hyperbolic space.  Geometriae Dedicata. {\bf 106}, 51-78 (2004)


\bibitem{kimp}  Kim, I., Parker, J. R.: Geometry of quaternionic hyperbolic
manifolds,  Math. Proc. Cambridge Philos. Soc. {\bf 135}, 291-320
(2003)


\bibitem{km} Kor\'anyi, A., Reimann, H. M.: The complex cross-ratio on the
Heisenberg group.  Enseign. Math. {\bf 33}, 291-300 (1987)



\bibitem{kel1}  Kellerhals, R.: Collars in PSL(2, H). Ann. Acad. Sci. Fenn.
Math. {\bf 26}, 51-72 (2001)

\bibitem{kel2} Kellerhals, R.: Quaternions and some global properties of 
hyperbolic 5-manifolds. Canad. J.
Math. {\bf 55}, 1080-1099 (2003)

\bibitem{mar} Markham, S.: Hypercomplex Hyperbolic Geometry.  
PhD Thesis, Univ. Durham, (2003)

\bibitem{mapa} Markham, S., Parker, J. R.:   Collars in complex and quaternionic
hyperbolic manifolds.  Geometriae Dedicata. {\bf 97}, 199-213 (2003)

\bibitem{mapa2} Markham, S., Parker, J. R.:  J\o rgensen's inequality for
metric spaces with applications to the octonions.
Advances in Geometry {\bf 7}, 19--38 (2007).


\bibitem{mey} Meyerhoff, R.: A lower bound for
the volume of hyperbolic 3-manifolds, Canad. J. Maths. {\bf 39},
1038-1056 (1987)

\bibitem{mos} Mostow, G.D.: Strong Rigidity of Locally Symmetric Spaces,
Annals of Maths.\ Studies {\bf 78}, Princeton 1973.

\bibitem{parsbt} Parker, J.R.: On the stable basin theorem, 
Canadian Mathematical Bulletin {\bf 47} 439-444 (2004). 

\end{thebibliography}
\end{document}